\documentclass{amsart}
\usepackage[mathscr]{eucal}
\usepackage{amssymb, amsmath,array, amscd}
\usepackage[OT2,T1]{fontenc}

\newtheorem{thm}{Theorem}
\newtheorem{thma}{Theorem}[section]

\newtheorem{lemma}{Lemma}[section]
\newtheorem{claim}{Claim}[section]

\newtheorem{definition}{Definition}

\newtheorem{rem}{Remark}[section]

\newtheorem{cor}{Corollary}[section]

\newtheorem{Cor}{Corollary}

\newtheorem{prob}{Problem}

\newtheorem{ex}{Example}

\def\wh{\widehat}

\def\de{\mathcal{D}}

\def\td{\tilde}

\def\Ga{\Gamma}

\def\p{\mathbb{P}}

\def\Te{\Theta}
\def\te{\theta}
\def\lim{\underset{z\to 0}{lim}\,}
\def\pa{\partial}
\def\sup{\supset}
\def\sub{\subset}
\def\De{\Delta}

\def\d{\delta}

\def\emp{\emptyset}
\def\ov{\overline}
\def\om{\omega}
\def\Om{\Omega}
\def\fa{\mathcal{F}}
\def\a{\alpha}
\def\be{\beta}
\def\R{\mathbb{R}}
\def\ep{\epsilon}
\def\C{\mathbb{C}}

\def\{{\lbrace}
\def\}{\rbrace}
\def\D{\mathbb{D}}

\def\la{\lambda}

\def\U{\mathcal{U}}

\def\g{\gamma}
\def\P{\mathcal{P}}
\def\N{\mathbb{N}}
\def\Si{\Sigma}
\def\te{\theta}
\def\si{\sigma}
\def\G{\mathcal G}

\def\Z{\mathbb{Z}}
\def\O{\mathcal{O}}
\def\*{\star}

\def\La{\Lambda}

\def\Q{\mathbb{Q}}

\def\var{\varphi}

\def\X{\mathcal{X}}

\begin{document}
\title[Transversal product singular set]{Local transversely product singularities}


\author[A. Lins Neto]{A. Lins Neto}
\address{IMPA, Est. D. Castorina, 110, 22460-320, Rio de Janeiro, RJ, Brazil}

\email{alcides@impa.br}

\thanks{}



\begin{abstract}
In the main result of this paper we prove that a codimension one foliation of $\p^n$, which is locally a product near every point of some codimension two component of the singular set, has a Kupka component. In particular, we obtain a generalization of a known result of Calvo Andrade and Brunella about foliations with a Kupka component.
\end{abstract}

\maketitle

\tableofcontents

\section{Basic definitions and results}

It is known that a holomorphic codimension one foliation on $\p^n$, $n\ge3$, with a Kupka component in the singular set has a rational first integral, which in homogeneous coordinates is of the form $P^k/Q^\ell$, where $P$ and $Q$ are generic homogeneous polynomials with $\ell.\,deg(P)=k.\,deg(Q)$.
The Kupka component for this specific example is the set $\Pi(P=Q=0)$, where $\Pi\colon\C^{n+1}\setminus\{0\}\to\p^n$ is the canonical projection (cf. \cite{br}, \cite{ca1} and \cite{ca2}).
The aim of this paper is to generalize this result for codimension one foliations with a {\it local transversely product component} in the singular set.
We will define this concept in a more general situation.

Let $\fa$ be a holomorphic foliation of dimension $k\ge2$ on a complex manifold $M$ of dimension $n\ge k+1$, with singular set $Sing(\fa)$. We say that $\fa$ is {\it a transversely product} at a point
$p\in Sing(\fa)$ if the germ $\fa_p$ of $\fa$ at $p$ is holomorphically equivalent to a product of a germ of singular foliation of dimension one with an isolated singularity by a regular foliation of dimension $k-1$. In other words, we can say that there exists a germ of submersion $\var\colon(M,p)\to(\C^{n-k+1},0)$ and a germ of an one dimensional foliation $\G$ at $0\in\C^{n-k+1}$, with $Sing(\G)=\{0\}$, such that $\fa_p=\var^*(\G)$.
In particular, the germ of the singular set of $\fa$ at $p$ is smooth of dimension $k-1$: $Sing(\fa_p)=\var^{-1}(0)$.

\begin{definition}\label{def:1}
{\rm We say that $\Ga$ is a {\it local transversely product component} (briefly l.t.p component) of $Sing(\fa)$ if $\Ga$ is an irreducible component of $Sing(\fa)$ and $\fa$ is a transversely product at all points of $\Ga$.}
\end{definition}

\begin{rem}\label{r:11}
{\rm If $\Ga$ is a l.t.p. component of $Sing(\fa)$ then it follows from the definition that:
\begin{itemize}
\item[a.] $\Ga$ is smooth. Let $dim_\C(\Ga)=m$.
\item[b.] There exists a singular one dimensional foliation $\G$, on a polydisc $V$ of $\C^{n-m}$, with an isolated singularity at $0\in V$, such that for any $p\in\Ga$ there exists a local chart $(U,z)$ around $p\in U$ satisfying the following conditions:
\begin{itemize}
\item[b.1.] $z=(x,y)\colon U\to\C^{n-m}\times \C^{m}$ with $x(U)=V$.
\item[b.2.] $\fa|_U=x^*(\G)$.
\end{itemize}
\end{itemize} 
In the chart $z=(x,y)$ the submersion of the definition is $\var=x\colon U\to V$ and the leaves of the non-singular foliation are the levels $x^{-1}(a)$, $a\in V$. Moreover, $\Ga\cap U=x^{-1}(0)$.

The germ of $\G$ at $0\in Q$ is called the {\it normal type} of $\fa$ along $\Ga$.
Remark that, if $T$ is a germ at $p\in\Ga$ of $n-m$ manifold transverse to $\Ga$ then the restricted foliation $\fa|_T$ is holomorphically equivalent to the normal type of $\fa$ along $\Ga$.

Moreover, since $\G$ is one dimensional we can assume that it is defined by a holomorphic vector field $X=\sum_{j=1}^{n-m}A_j(x)\frac{\pa}{\pa x_j}$, or by the $(n-m-1)$-form
\begin{equation}\label{eq:1}
\eta=i_Xdx_1\wedge...\wedge dx_{n-m}=\sum_{j=1}^{n-m}(-1)^{j-1}A_j(x)\,dx_1\wedge...\wedge\wh{dx_j}\wedge...\wedge dx_{n-m}\,,
\end{equation}
where in (\ref{eq:1}) $\wh{dx_j}$ means omission of $dx_j$ in the product.
The form $\eta$, considered as a form on $U$ in the coordinates $(x,y)$, defines $\fa|_U$.

Let us see some examples:}
\end{rem}

\begin{ex}\label{ex:1}
{\rm Recall that a point $p\in M$ is a Kupka singularity of the foliation $\fa$ if $p\in Sing(\fa)$ and $\fa$ is represented in a neighborhood of $p$ by an integrable $(n-k)$-form $\eta$ such that $d\eta(p)\ne0$. The form $d\eta$ defines a $k+1$ distribution $\de=ker(d\eta)$ in a neighborhood of $p$, where
\[
\de(q)=ker(d\eta(q)):=\{v\in T_qM\,|\,i_v(d\eta(q))=0\}\,.
\]
The distribution $\de$ is integrable and defines a regular foliation of dimension $k-1$ in a neighborhood of $p$. There exists a local chart $(U,z=(x,y))$, where $x=(x_1,...,x_{n-k+1})\colon U\to \C^{n-k+1}$, $y\colon U\to\C^{k-1}$ and $z(p)=(0,0)$, such that $d\eta=dx_1\wedge...\wedge dx_{n-k+1}$. In this case, the form $\eta$ can be written as $\eta=i_Xdx_1\wedge...\wedge dx_{n-k+1}$, where
\[
X=\sum_{j=1}^{n-k+1}A_j(x)\,\frac{\pa}{\pa x_j}
\]
defines the normal type of $\fa_p$.

Note that $d\eta=\De(X)\,dx_1\wedge...\wedge dx_{n-k+1}$, where
\[
\De(X)=div(X)=\sum_{j=1}^{n-k+1}\frac{\pa A_j}{\pa x_j}\,,
\]
so that $\De(X)\equiv1$.}
\end{ex}

\begin{definition}\label{def:2}
{\rm We say that $K$ is a Kupka component of a foliation $\fa$ (of dimension $k\ge2$) if $K$ is a l.t.p. component of $Sing(\fa)$ and the normal type of $\fa$ along $K$ is of Kupka type.}
\end{definition}

\begin{ex}\label{ex:2}
{\rm Let $P$ and $Q$ be homogeneous polynomials on $\C^{n+1}$, $n\ge3$, where $def(P)=p$ and $deg(Q)=q$. 
The levels of the rational function $f=\frac{P^q}{Q^p}$ define a singular foliation of $\p^n$, that will be denoted by $\fa(P,Q)$.
We say that $P$ and $Q$ are transverse if the set
\[
\{z\in \C^{n+1}\,|\,P(z)=Q(z)=0\,\text{and}\,\,dP(z)\wedge dQ(z)=0\}
\]
is either $\{0\}$, or empty (if $p=q=1$).
If $P$ and $Q$ are transverse then the subset $\Ga$ of $\p^n$ defined in homogeneous coordinates by $(P=Q=0)$ is a Kupka component on $\fa(P,Q)$.
The normal type of $\fa(P,Q)$ at the points of $\Ga$ is given by the linear vector field $X=p.\,x\frac{\pa}{\pa x}+q.\,y\frac{\pa}{\pa y}$.

In fact, the following result is known (cf. \cite{br}, \cite{ca1}, \cite{ca2} and \cite{lc}):

\begin{thma}\label{t:11}
Let $\fa$ be a holomorphic foliation of codimension one on $\p^n$, $n\ge3$. If $\fa$ has a Kupka component then $\fa=\fa(P,Q)$, where $P$ and $Q$ are transverse polynomials.
\end{thma}}
\end{ex}

\begin{ex}\label{ex:3}
{\rm Example \ref{ex:2} admits the following generalization: let $P_1,...,P_m$ be homogeneous polynomials on $\C^{n+1}$ with $deg(P_j)=d_j$, $1\le j\le m$. Assume that
$n\ge m+1\ge4$ and that $P_1,...,P_m$ are transverse; i.e. the set
\[
\{z\in\C^{n+1}\,|\,P_1(z)=...=P_m(z)=0\,\text{and}\,dP_1(z)\wedge...\wedge dP_m(z)=0\}
\]
is either $\{0\}$, or empty (if $d_1=...=d_m=1$).
Let $(k_1,...,k_m)\in\N^m$ be such that $gcd(k_1,...,k_m)=1$ and $k_1.\,d_1=...=k_m.\,d_m$. The levels of the rational map $\P\colon\p^n\to\p^{m-1}$, defined by
\[
\P:=[P_1^{k_1}:...:P_m^{k_m}]\,,
\]
define a foliation of codimension $m-1$ on $\p^n$, denoted by $\fa(P_1,...,P_m)$. If $P_1,...,P_m$ are transverse then the set $\Ga\sub\p^n$, defined in homogeneous coordinates by $(P_1=...=P_m=0)$, is a Kupka component of $\fa(P_1,...,P_m)$. The normal type of $\fa(P_1,...,P_m)$ at the poins of $\Ga$ is given by the linear vector field
$S=\sum_{j=1}^m d_j\,x_j\frac{\pa}{\pa x_j}$\,.

A natural problem is the following:
\begin{prob}
{\rm Let $\fa$ be a holomorphic codimension $m-1$ on $\p^n$, where $n\ge m+1\ge4$. Assume that $\fa$ has a Kupka component $\Ga$. Are there transverse homogeneous polynomials $P_1,...,P_m$ on $\C^{n+1}$ such that $\fa=\fa(P_1,...,P_m)$ and $\Ga$ is defined by $(P_1=...=P_m=0)$?}
\end{prob}
Some partial results about this problem were proved (see for instance \cite{ca3} and \cite{moa}).}
\end{ex}

In this paper we generalize theorem \ref{t:11}:

\begin{thm}\label{t:1}
Let $\fa$ be a holomorphic foliation of codimension one on $\p^n$, $n\ge3$. Assume that $\fa$ has a l.t.p. component $\Ga$. Then $\Ga$ is a Kupka component of $\fa$. In particular, $\fa$ is like in example \ref{ex:2}.
\end{thm}

Let us state some consequences of theorem \ref{t:1}.

\begin{Cor}\label{c:1}
Let $\fa$ be a codimension one holomorphic foliation on $\p^n$, $n\ge4$. Assume that there is a linear embedding $i\colon\p^3\to\p^n$ such that $i^*(\fa)$ has a l.t.p component.
Then $\fa$ has a rational first integral that can be written in homogeneous coordinates as $P^q/Q^p$, where $P$ and $Q$ are homogeneous polynomials on $\C^{n+1}$ with $deg(P)=p$ and $deg(Q)=q$.
\end{Cor}

The proof of corollary \ref{c:1} is based in the fact that if there exists a linear embedding $i\colon\p^3\to\p^n$ such that $i^*(\fa)$ has a first integral then $\fa$ has also a first integral (see \cite{lc3}).

\begin{Cor}\label{c:2}
Let $\fa$ be a codimension one foliation on $\p^n$, $n\ge3$. Assume that all components of its singular set are l.t.p. Then $\fa$ has degree zero: the first integral of corollary \ref{c:1} is of the form $L_2/L_1$, where $L_1$ and $L_2$ are linear.
\end{Cor}

\begin{Cor}\label{c:3}
Let $\eta$ be an integrable 2-form on $\C^n$, $n\ge4$, with homogeneous coefficients of the same degree $d\ge1$. Then $dim_\C(sing(\eta))\ge1$.
\end{Cor}

\begin{rem}
{\rm Corollary \ref{c:3} was proved in \cite{lc1} in the case $n=4$. We would like to observe that the assertion is not true in the case of distributions of $\C^4$.
The following example, due to Krishanu and Nagaraj \cite{KN}: define a 2-form $\te$ on $\C^4$ by
\[
\te=x_3^2\,dx_2\wedge dx_3-x_1^2\,dx_3\wedge dx_1+(x_1\,x_2+x_3\,x_4)\,dx_1\wedge dx_2+
\]
\[
[x_4^2\,dx_1+x_2^2\,dx_2+(x_1\,x_2-x_3\,x_4)\,dx_3]\wedge dx_4
\]
has $Sing(\te)=\{0\}$ and satisfies $\te\wedge\te=0$. Hence, it generates a distribution of codimension two on $\C^4\setminus\{0\}$. This distribution is not integrable.
}
\end{rem}

Theorem \ref{t:1} motivates the following problem:

\begin{prob}
{\rm Let $\fa$ be a holomorphic foliation on $\p^n$ of codimension $\ge2$ and dimension $\ge2$. Assume that $\fa$ has a l.t.p. component $\Ga$. Is $\Ga$ a Kupka component of $\fa$?}
\end{prob}

A crucial point of our proof of theorem \ref{t:1} is the Camacho-Sad theorem on the existence of a separatrix for germs of holomorphic vector fields on $(\C^2,0)$ \cite{cs}. The same type of argument cannot be used in the general case: there are examples of germs of vector fields on $(\C^m,0)$, $m\ge3$, without separatrices \cite{gl}.

\vskip.1in
 
The proof of theorem \ref{t:1} will be done in section \ref{ss:2}. Since this proof is technical, in section \ref{ss:21} we give an idea of the proof by stating the main objects and results that will be used. In sections \ref{ss:22} and \ref{ss:23} we will prove the main auxiliary results used in the proof and stated in section \ref{ss:21}.
Section \ref{ss:3} is dedicated to the proof of corollaries \ref{c:2} and \ref{c:3}.

\vskip.1in

{\bf Aknowledgement.} I would like to thank J. Vitório Pereira for helpfull conversations that suggested me a simplification of the proof of lemma \ref{l:23}. 

\section{Proof of theorem \ref{t:1}}\label{ss:2}

\subsection{Preliminaries and idea of the proof}\label{ss:21}
Let $\fa$ be a codimension one foliation on $\p^n$, $n\ge3$, with a l.t.p. component $\Ga\sub Sing(\fa)$.
The definition implies that $cod_\C(\Ga)=2$, so that the transversal type of $\fa$ at the poins of $\Ga$ is a germ of singular foliation at $(\C^2,0)$ with an isolated singularity at $0\in\C^2$ (see remark \ref{r:11} and example \ref{ex:1}).
We can assume that this transversal type is given by germ at $0\in\C^2$ of vector field $X=X_1(x,y)\frac{\pa}{\pa x}+X_2(x,y)\frac{\pa}{\pa y}$, where $X_1$, $X_2\in\O_2$ and $X_1(0,0)=X_2(0,0)=0$.
Recall that $\Ga$ is a Kupka component if, and only if, we have $Tr(DX(0))\ne0$, where
\[
Tr(DX(0)):=\frac{\pa X_1}{\pa x}(0)+\frac{\pa X_2}{\pa y}(0)
\]
is the trace of the linear part $DX(0)$ of $X$ at $0\in\C^2$.
In this case, as we have pointed out before, $\fa$ is like in example \ref{ex:2} (see theorem \ref{t:11}).

Another useful ingredient is the normal Baum-Bott index of the component $\Ga$, that we will denote as $BB(\fa,\Ga)$. Since $\Ga$ is a l.t.p. component of $Sing(\fa)$ then $BB(\fa,\Ga)$ coincides with the Baum-Bott of $X$ at the singularity $0$ of $X$, denoted by $BB(X,0)$ (see \cite{lc2}) (for the definition of $BB(X,0)$ see \cite{br1}). In lemma 3.4 of \S 3.2 of \cite{lc2} it is proven that if $BB(\fa,\Ga)\ne0$ and $DX(0)\not\equiv0$ then $\Ga$ is a Kupka component and we are done.

One of the tools used in the proof of lemma 3.4 of \cite{lc2} is the existence of a smooth analytic separatrix along $\Ga$.
Below we define the concept of separatrix in a way that will be used in the proof of theorem \ref{t:1}.

\begin{definition}\label{def:3}
{\rm Let $\fa$ be a holomorphic foliation of dimension $k$ on a n dimensional compact complex manifold, $2\le k<n$, and $\Ga$ be l.t.p. component of $\fa$ (recall that $dim(\Ga)=k-1$). A {\it separatrix of dimension $\ell$ along} $\Ga$ of $\fa$, where $k\le\ell<n$, is a germ of $\ell$ analytic manifold along $\Ga$ which is $\fa$-invariant in the sense that:
\begin{itemize}
\item[a.] $\Si\sup\Ga$.
\item[b.] $\Si\setminus\Ga$ is contained in an union of leaves of $\fa$.
\end{itemize}}
\end{definition}

\begin{rem}\label{r:21}
{\rm Let $\Si$ be a separatrix of $\fa$ of dimension $\ell$ along $\Ga$, as in definition \ref{def:3}. Fix $p\in\Ga$ and $(x,y)\colon U\to\C^{n-k+1}\times\C^{k-1}$, a local coordinate system around $p$ as in remark \ref{r:11}. It follows from the definition that $x^{-1}(x(\Si\cap U))=\Si\cap U$.

Let $T$ is a germ at $p$ of a $n-k+1$ dimensional manifold transverse to $\Ga$. As we have observed before, $\fa|_T$ is equivalent to the normal type of $\fa$ along $\Ga$. In particular, the intersection $\Si\cap T$ is invariant by $\fa|_T$.

In the case of theorem \ref{t:1}, where $\fa$ has codimension one, then $dim(T)=2$ and the normal type is a germ $\G$ of one dimensional foliation on $(\C^2,0)$. In this case $\Si\cap T$ is a finite number of analytic separatrixes of $\G$ as considered in \cite{cs}. The next result will be used in proof of theorem \ref{t:1}.}
\end{rem}

\begin{lemma}\label{l:21}
Let $\fa$ be a holomorphic codimension one foliation on a compact complex manifold $M$, where $dim(M)\ge3$, and $\Ga$ be a l.t.p. component of $Sing(\fa)$. If the normal type of $\fa$ along $\Ga$ is not equivalent to the radial foliation of $(\C^2,0)$ then $\fa$ admits an irreducible separatrix $\Si$ along $\Ga$ with $dim(\Si)=n-1$.  
\end{lemma}

Recall that the radial foliation of $\C^2$ is defined by the form $x\,dy-y\,dx$ and their leaves are the straight lines through $0$.
Lemma \ref{l:21} will be proved in section \ref{ss:22}.

\vskip.1in

From now on, in this section, we will assume that $\fa$ is a codimension one holomorphic foliation on the compact manifold $M$, $dim_\C(M)\ge3$, with a l.t.p. component $\Ga$ and with a separatrix $\Si$ along $\Ga$, $dim(\Si)=n-1$.
Next we will introduce the {\it normal bundle of $\Si$ along $\Ga$}.

Since $dim(\Si)=n-1$ we can find a Leray covering $\U=(U_\a)_{\a\in A}$ of $\Ga$ by open sets and two collections $f=(f_\a)_{\a\in A}$ and $g=(g_{\a\be})_{U_\a\cap U_\be\ne\emp}$ with the following properties:
\begin{itemize}
\item[(a).] $f_\a\in\O(U_\a)$, $\forall\a\in A$, and $f_a=0$ is a reduced equation of $\Si\cap U_\a$.
\item[(b).] $g_{\a\be}\in\O^*(U_\a\cap U_\be)$ and $f_\a=g_{\a\be}.\,f_\be$ on $U_\a\cap U_\be\ne\emp$.
\end{itemize}

Of course $g=(g_{\a\be})_{U_\a\cap U_\be\ne\emp}$ is a multiplicative cocycle.
We define the normal bundle of $\Si$ along $\Ga$ as the line bundle on $Picc(\Ga)$ induced on a tubular neighborhood $U\sub\bigcup_\a U_\a$ by the cocycle $g=(g_{\a\be})_{U_\a\cap U_\be\ne\emp}$.

It will be denoted by $N_\Si$. Let $C_1(N_\Si)$ be the first Chern class of $N_\Si$, considered as an element of $H^2(U,\R)$ via the homomorphism
$H^2(U,\Z)\to H^2(U,\R)\simeq H^2_{Dr}(U)$ induced by the inclusion $\Z\to\R$.
 
As we have seen in remark \ref{r:11}, the normal type of $\fa$ along $\Ga$ can be represented by a germ at $0\in\C^2$ of holomorphic vector field
$X=A_1(x,y)\,\frac{\pa}{\pa x}+A_2(x,y)\,\frac{\pa}{\pa y}$ with an isolated at $0$. When we intersect $\Si$ with a germ of transversal section $T\simeq(\C^2,0)$ we obtain a separatrix of $X$, say $\g:=\Si\cap T$ (in general $\g$ is not irreducible).
Let $f\in\O_2$ be a reduced analytic equation of $\g$. Since $\g$ is $X$-invariant we can write
\begin{equation}\label{eq:}
X(f)=h.\,f\,\,,\,\,\text{where}\,\,h\in\O_2\,.
\end{equation}

\begin{lemma}\label{l:22}
In the above situation, if $h(0)=0$ then $C_1(N_\Si)=0$.
\end{lemma}

On the other hand, we have the following:

\begin{lemma}\label{l:23}
If the ambient space is $M=\p^n$, $n\ge3$, then $C_1(N_\Si)\ne0$.
In particular, if $X(f)=h.\,f$ then $h(0)\ne0$.
\end{lemma}

As a consequence of lemma \ref{l:23} we get the following:

\begin{cor}\label{c:21}
If $M=\p^n$, $n\ge3$, then $\Si$ is a Kupka component of $\fa$.
\end{cor}

In particular, theorem \ref{t:11} will imply theorem \ref{t:1}.
Lemma \ref{l:21} will be proved in the next section.

\vskip.1in

\subsection{Proof of lemma \ref{l:21}.}\label{ss:22}

Let $\fa$ be a holomorphic codimension one foliation on a compact complex manifold $M$ with $dim(M)\ge3$. Assume that $\fa$ has a l.t.p. component $\Ga$ with normal type $\G$, where $\G$ is a germ of foliation on $(\C^2,0)$ with an isolated singularity at $0\in\C^2$.
As before, we will assume that $\G$ is the foliation defined by a germ at $(\C^2,0)$ of vector field $X=X_1\frac{\pa}{\pa x}+X_2\frac{\pa}{\pa y}$ with an isolated singularity at the origin of $\C^2$.
The germ of foliation $\G$ can be defined also by the 1-form 
\[
\om=i_X(dx\wedge dy)=X_1\,dy-X_2\,dx\,,
\]
so that, $d\om(0)=Tr(DX(0))\,dx\wedge dy$. We can assume that $\om$ has a representative, denoted by $\td\om$, defined in the polydisc $Q=\D^2$ with an isolated singularity at $0\in\D^2$.

By the definition of l.t.p. component, we can find a covering $\U=(U_\a)_{\a\in A}$ of $\La$ by open sets biholomorphic to polydiscs, a collection of local charts $\left((z_\a,U_\a)\right)_{\a\in A}$ and a multiplicative cocicle $(k_{\a\be})_{U_\a\cap U_\be\ne\emp}$ with the following properties:
\begin{itemize}
\item[1.] $z_\a=(x_\a,y_\a)\colon U_\a\to\C^2\times\C^{n-2}$, where $x_\a(U_\a)=Q$ and $\Ga\cap U_\a=x_\a^{-1}(0)$, $\forall\a\in A$.
\item[2.] $\fa|_{U_\a}$ is defined by the integrable 1-form $\td\om_\a:=x_\a^*(\td\om)$. The germ of $\td\om_\a$ along $\Ga\cap U_\a$ will be denoted by $\om_\a$.
\item[3.] $\td\om_a=k_{\a\be}.\,\td\om_\be$ on $U_\a\cap U_\be\ne\emp$.
\end{itemize}
We will assume that $\U$ satisfies the following:
\begin{itemize}
\item[4.] If $U_\a\cap U_\be\ne\emp$ then $\Ga\cap U_\a\cap U_\be\ne\emp$ and connected. 
\end{itemize}

\begin{rem}\label{r:22}
{\rm Given $\a,\be\in A$ such that $\Ga\cap U_\a\cap U_\be\ne\emp$ we can construct a germ $f_{\a\be}\in Diff(\C^2,0)$ as follows: fix $p\in\Ga\cap U_\a\cap U_\be$ and a germ of plane $T=T_{\a,\be}\simeq(\C^2,p)$ transverse to $\Ga$ at $p$. Note that $x_\a|_T,x_\be|_T\colon (T,p)\to(\C^2,0)$ are biholomorphisms. Therefore, we define
\[
f_{\a\be}=x_\a\circ(x_\be|T)^{-1}=x_\a|T\circ(x_\be|T)^{-1}\in Diff(\C^2,0)\,.
\]
Since $\om_\a|_T=(x_\a|T)^*(\om)$, $\om_\be|_T=(x_\be|_T)^*(\om)$ and $\om_a=k_{\a\be}.\,\om_\be$ we get $f_{\a\be}^*(\om)=h_{\a\be}.\,\om$, where
\[
h_{\a\be}=k_{\a\be}|_T\circ(x_\be|_T)^{-1}\in\O_2^*\,.
\]
The biholomorphism $f_{\a\be}$ can be interpreted as the {\it glueing map} of $\fa|_{U_\be}$ with $\fa|_{U_\be}$. }
\end{rem}

From now on, we fix a collection of germs $(f_{\a\be})_{U_\a\cap U_\be\ne\emp}$ as above.

\begin{lemma}\label{l:24}
$\fa$ admits a separatrix $\Si$ along $\Ga$ if, and only if, $X$ (or $\om$) has a separatrix $\g$ (not necessarily irreducible) such that $f_{\a\be}(\g)=\g$ for all
$\Ga\cap U_\a\cap U_\be\ne\emp$.
\end{lemma}

{\bf Notation.}
We will say that the separatrix $\g$ of $\G$ generates the separatrix $\Si$ of $\fa$.

{\it Proof.} Assume that $X$ has a separatrix $\g$ such that $f_{\a\be}(\g)=\g$ for all $\Ga\cap U_\a\cap U_\be\ne\emp$.
Given $\a\in A$ define $\Si_\a:=x_\a^{-1}(\g)$. We assert that if $U_\a\cap U_\be\ne\emp$ then $\Si_\a\cap U_\be=\Si_\be\cap U_\a$.
In fact, let $(x_\a,y_\a)$, $(x_\be,y_\be)$ and $T$ be as before. Then
\[
\Si_\a\cap T=x_\a^{-1}(\g)\cap T=(x_\a|_T)^{-1}(\g)=(x_\a|_T)^{-1}(f_{\a\be}(\g))=(x_\be|_T)^{-1}(\g)=\Si_\be\cap T\,.
\]
This, of course, implies the assertion.
In particular, the local separatrices $\Si_a$ glue together forming a global separatrix $\Si$ along $\Ga$ such that $\Si\cap U_\a=\Si_\a$, $\forall\a\in A$.

We leave the converse to the reader.
\qed

\begin{definition}\label{def:4}
{\rm Let $\G$ be a germ of foliation at $(\C^2,0)$ with an isolated singularity at $0$. We say that a separatrix $\g$ of $\G$ is {\it distinguished} if for any $f\in Diff(\C^2,0)$ such that $f^*(\G)=\G$ then $f(\g)=\g$.}
\end{definition}

\begin{lemma}\label{l:25}
Let $\G$ be a germ of foliation at $(\C^2,0)$ with an isolated singularity at $0$ which is not equivalent to the radial foliation.
Then $\G$ has a distinguished separatrix.  
\end{lemma}

{\it Proof.}
In the proof we use Seidenberg's resolution theorem \cite{sd}.
Let $S$ be a smooth complex surface and $\G$ be a foliation by curves on $S$. Given $p\in Sing(\G)\sub S$ we denote $Diff(S,p)$ the set of germs at $p\in S$ of biholomorphisms $f\colon(S,p)\to S$ with a fixed point at $p$. Assume that the germ of $\G$ at $p$ is defined by a germ of holomorphic vector field $X$ with an isolated singularity at $p$. We use also the notations
\[
Diff_{\G}(S,p)=\{f\in Diff(S,p)\,|\,f^*(\G)=\G\}\,.
\]
\begin{rem}\label{r:23}
{\rm
Note that:
\begin{itemize}
\item[(1).] Given $f\in Diff_\G(S,p)$, then $f^*(X)=h_X.\,X$, where $h_X\in\O^*_p$.
\item[(2).] $Diff_\G(S,p)$ is a sub-group of $Diff((S,p)$.
\item[(3).] Given $f\in Diff_\G(S,p)$ and an irreducible separatrix $\g$ of $\G$ through $p$ then $f(\g)$ is also a separatrix of $\G$ through $p$.
\end{itemize}}
\end{rem}

Let $Sep(\G)$ be the set of irreducible separatrices of $\G$ through $p$.
By (3) of remark \ref{r:23}, $Diff_\G(S,p)$ acts in $Sep(\G)$ as $(f,\d)\in Diff_\G(S,p)\times Sep(\G)\to f(\d)\in Sep(\G)$.
The idea of the proof is to find a finite subset $G_o:=\{\g_1,...,\g_k\}\sub Sep(\G)$ such that $f(G_o)\sub G_o$ for all $f\in Diff_\G(S,p)$.
In this case, the set $\g:=\{f(\g_1)\,|\,f\in Diff_\G(S,p)\}\sub G_o$ contains finitely many irreducible separatrices of $\G$ through $p$ and can be considered as a germ of curve through $p$ such that $f(\g)=\g$ for all $f\in Diff_\G(S,p)$, and so $\g$ is a distinguished separatrix of $\G$ through $p$.
Let us prove the existence of the finite set $G_o$.

\vskip.1in

First of all, we observe that there are two possibilities for the foliation $\G$:
\begin{itemize}
\item[I.] $\G$ has finitely many irreducible separatrices through $p$.
This case is trivial and the details are left to the reader.
\item[II.] $\G$ has infinitely many irreducible separatrices through $p$. Let us prove lemma \ref{l:25} in this case.
\end{itemize}

We will consider a blowing-up process used to resolve the foliation $\G$ (see \cite{cs}).
The first case, is when $\G$ has a simple singularity at $p$ and no blowing-ups are needed in the process.
Let $\la_1$ and $\la_2$ be the eigenvalues of $DX(p)$.
The singularity is simple if:
\begin{itemize}
\item[(a).] $\la_1.\,\la_2\ne0$ and $\frac{\la_2}{\la_1}\notin\Q_+$.
\item[(b).] $\la_1\ne0$ and $\la_2=0$ (or vice-versa). In this case, $p$ is a saddle-node. 
\end{itemize}
In both cases $\G$ has one or two separatrices through $p$ and so lemma \ref{l:25} is true.

When the singularity is not simple,
Seidenberg's theorem says that after a finite process of blowing-ups $\Pi\colon(\td{S},E)\to (S,p)$ then all the singularities of the strict transform $\Pi^*(\G)$ in the exceptional divisor $E$ are simple.
The blowing-up process $\Pi$ can be considered as a composition of pontual blowing-ups
\begin{equation}\label{eq:3}
(\td{S},E):=(\td{S}_k,E_k)\overset{\Pi_k}\longrightarrow(\td{S}_{k-1},E_{k-1})\overset{\Pi_{k-1}}\longrightarrow...\overset{\Pi_{2}}\longrightarrow(\td{S}_1,E_1)
\overset{\Pi_{1}}\longrightarrow(\td{S}_0,E_0)=(S,p)
\end{equation}
where in the $j^{th}$ step $\Pi_j\colon(\td{S}_j,E_j)\to(\td{S}_{j-1},E_{j-1})$, $j\ge2$, we blow-up in a point $p_{j-1}\in E_{j-1}$. The exceptional divisor obtained in this step will be denoted as $\p^1\simeq\td{E}_j\sub E_j$, so that $\Pi_j(\td{E}_j)=p_{j-1}$. We use also the notation $\td\Pi_j:=\Pi_1\circ...\circ\Pi_j$. 
We will denote also $\td\G_j:=\td\Pi^*(\G)$.
The point $p_{j-1}\in E_{j-1}$ is chosen between the non simple singularities of $\td{\G}_{j-1}$ on $E_{j-1}$.
Seidenberg's theorem can be stated as follows
\begin{thma}
It is possible to choose a blowing-up process as a above in such a way that all singularities of the strict transform $\td\G_k=\td\Pi_k^*(\G)$ are simple. 
\end{thma}

\begin{rem}
{\rm 
There are two possibilities in each step $\Pi_j\colon(\td{S}_j,\td{E}_j)\to(\td{S}_{j-1},p_{j-1})$.
We assume that $p_{j-1}$ is a non simple singularity of $\td\G_{j-1}$.
Let $X_{j-1}$ be a germ at $p_{j-1}$ of holomorphic vector field that represents the germ of $\td\G_{j-1}$ at $p_{j-1}$. Let
$X_\nu=P_\nu(x,y)\frac{\pa}{\pa x}+Q_\nu(x,y)\frac{\pa}{\pa y}$ be the first non-zero jet of $X_{j-1}$ at $p_{j-1}$, where $P_\nu$ and $Q_\nu$ are homogeneous polynomials of degree $\nu\ge1$. Set $F_{\nu+1}(x,y)=x.\,Q_\nu(x,y)-y\,P_\nu(x,y)$.
\begin{itemize}
\item[(i).] If $F_{\nu+1}\not\equiv0$ then $F_{\nu+1}$ is homogeneous of degree $\nu+1$ and the blowing-up is called {\it non-dicritical}. The divisor $\td{E}_j$ is invariant for the foliation $\td\G_j$ and the singularities of $\td\G_j$ on $\td{E}_j$ are the directions correspondent to the directions defined by $F_{\nu+1}(x,y)=0$.
\item[(ii).] If $F_{\nu+1}\equiv0$ then $X_\nu=F_{\nu-1}(x,y)\,R$, where $R=x\frac{\pa}{\pa x}+y\frac{\pa}{\pa y}$ is the radial vector field in $\C^2$ and $F_{\nu-1}$ is homogeneous of degree $\nu-1$. In this case, the blowing-up is called {\it dicritical}. The divisor $\td{E}_j$ is non-invariant for $\td\G_j$ and it is transverse to $\td{E}_j$ outside the set $V_j\sub\td{E}_j$ corresponding to the directions defined by the equation $F_{\nu-1}(x,y)=0$.

If $\nu=1$ then $\td\G_{j-1}$ is equivalent to the radial foliation at $p_{j-1}$. We will say that $p_{j-1}$ is a {\it radial} singularity of $\td\G_{j-1}$.
If $p_{j-1}$ is not radial for $\td\G_{j-1}$ then $V_j\ne\emp$ and we can divide it into two disjoint subsets $V_j=\tau_j\cup\si_j$, where
\begin{itemize}
\item{} $\si_j=Sing(\td\G_j)\cap\td{E}_j$.
\item{} $\tau_j=V_j\setminus Sing(\td\G_j)$. We call $\tau_j$ the set of tangencies of $\td\G_j$ with $\td{E}_j$. 
\end{itemize}
\end{itemize}
Remark also that $Sep(\G)$ is finite if, and only if, all blowing-ups in the process are non-dicritical.

Since in the blowing-up process, in each step, $1\le j\le k$, we blow-up in some non-simple singularity of $\td\G_{j-1}$, if $\td{E}_j$ is dicritical, at the end the tangencies $\tau_j$ "survive", in the sense that there exists a set $\tau\sub E_k$ such that for any $1\le j<k$ such that $\tau_j\ne\emp$ then
\[
\tau_j\sub\Pi_k\circ...\circ\Pi_{j+1}(\tau)
\]} 
\end{rem}
\vskip.1in
For each $1\le j\le k$ denote by $Diff(\td{S}_j,E_j)$ the set of germs of biholomorphisms $f\colon(\td{S}_j,E_j)\to(\td{S}_j,E_j)$.

\begin{definition}\label{def:5}
{\rm We say that $f\in Diff(S,p)$ can be lifted to $Diff(\td{S}_j,E_j)$ if there exists a germ of biholomorphism $\td{f}_j\in Diff(\td{S}_j,E_j)$ such that the diagram below commutes:
\[
\begin{matrix}
(\td{S}_j,E_j)&\overset{\td{f}_j}\longrightarrow&(\td{S}_j,E_j)\\
\td\Pi_j\downarrow&\,&\downarrow\td\Pi_j\\
(S,p)&\overset{f}\longrightarrow&(S,p)\\
\end{matrix}
\]}
\end{definition}

\begin{rem}
{\rm Observe that, if the lift $\td{f}_j$ of $f$ exists then it is unique. When $j=1$ (just one blowing-up) the lifting exists for any $f\in Diff(S,p)$, but if $j\ge2$ then there are germs $f\in Diff(S,p)$ that can not be lifted to $Diff(\td{S}_j,E_j)$. However, we have the following:} 
\end{rem}

\begin{claim}\label{cl:21}
The blowing-up process can be done in such a way that any $f\in Diff_\G(S,p)$ can be lifted to the last step in an unique $\td{f}=\td{f}_k\in Diff(\td{S}_k,E_k)$. Moreover, $\td{f}$ preserves $\td\G_k$ in the sense that $\td{f}^*(\td\G_k)=\td\G_k$.
\end{claim}

{\it Proof.}
We say that the $j^{th}$ step of the blowing-up process is admissible if any $f\in Diff_\G(S,p)$ has a lifting $\td{f}_j\in Diff(\td{S}_j,E_j)$.
We will obtain by induction a blowing-up process, as in (\ref{eq:3}), for which there are steps $1=\ell_1<\ell_2<...<\ell_r=k$ such that the $\ell_j^{th}$ step is admissible, for any $1\le j\le r$, and $\td\Pi_k\colon(\td{S}_k,E_k)\to(S,p)$ is a resolution of the folition $\G$.

First of all, the first step is admissible, because any $f\in Diff(S,p)$ admits a lifting $\td{f}_1\in Diff(\td{S}_1,E_1)$.

Assume that we have found some process for wich the $\ell:=\ell_s$ step is admissible, $\ell\ge1$, so that any $f\in Diff_\G(S,p)$ admits a lifting $\td{f}=\td{f}_\ell\in Diff(\td{S}_\ell,E_\ell)$ satisfying $\td{f}^*(\td\G_\ell)=\td\G_\ell$.
Given $f\in Diff_\G(S,p)$, with lifting $\td{f}$, and $q\in E_\ell$ then $\td{f}$ is an equivalence between the two germs of $\td\G_\ell$ at $q$ and at $\td{f}(q)$.
In particular, $\td{f}$ preserves the set of non simple singularities of $\td\G_\ell$.

If $\td\G_\ell$ is not a resolution of $\G$ then it has at least one non simple singularity $q_1$.
Let $Sat(q_1)=\{\td{f}(q_1)\,|\,f\in Diff_\G(S,p)\}=\{q_1,...,q_m\}$.
We then blow-up once at all points $q_j\in Sat(q_1)$, passing from the $\ell=\ell_s$ step to the $\ell_{s+1}:=\ell_s+m$ step directly.
Let $\wh{E}_j$ be the divisor obtained by the blowing-up at $q_j$.

Given $f\in Diff_\G(S,p)$ and its lifting $\td{f}_s$, let $\td{f}_s(q_j)=q_{i(j)}$, $1\le j\le m$. Then we can obtain a lifting $\td{f}_{s+1}$ of $\td{f}_s$ 
such that $\td{f}_{s+1}(\wh{E}_j)=\wh{E}_{i(j)}$, $1\le j\le m$.
By Seidenberg's theorem this process must end at some step, when the final foliation $\td\G_k=\td\Pi_k^*(\G)$ has all singularities simple.
\qed 

\vskip.1in

Let us finish the proof of lemma \ref{l:25}. We will consider two cases:
\begin{itemize}
\item[(1).] There is $q_1\in Sing(\td\G_k)\cap E_k$ that has some separatrix $\td\g$ not contained in $E_k$.
\item[(2).] All the separatrices of the singularities of $\td\G_k$ are contained in $E_k$.
\end{itemize}

In the first case, let $Sat(q_1)=\{\td{f}(q_1)\,|\,f\in Diff_\G(S,p)\}=\{q_1,...,q_m\}$. Given $f\in Diff_\G(S,p)$ and $q_j=\td{f}(q_1)$ then $\td{f}(\td\g):=\td\g_f$ is a separatrix of $\td\G_k$ not contained in $E_k$. Since $\td\g_f$ is not contained in $E_k$, its immage $\g_f:=\td\Pi_k(\td\g_f)$ is a separatrix of $\G$ through $p$. Moreover, since $q_1,...,q_m$ are simple singularities of $\td\G_k$ the set $\{\td\g_f\,|\,f\in Diff_\G(S,p)\}$ is finite.
Therefore, if we set
\begin{equation}\label{eq:4}
\g=\bigcup_{f\in Diff_\G(S,p)}\,\g_f
\end{equation}
then $f(\g)=\g$, $\forall f\in Diff_\G(S,p)$, and $\g$ is a distinguished separatrix of $\G$.

In the second case necessarily there are dicritical irreducible divisors of $\td\G_k$, say $\td{E}_1,...,\td{E}_m$, contained in $E_k$ (by Camacho-Sad theorem).
This case will be divided into two sub-cases:
\begin{itemize}
\item[(2.1).] The set of tangencies $\tau$ is not empty.
\item[(2.2).] $\tau=\emp$.
\end{itemize}

In case (2.1) let $q_o\in \tau$ and $\td\g_{\,id}$ be the leaf of $\td\G_k$ through $q_o$. Then, for any $f\in Diff_\G(S,p)$ we have $q_f:=\td{f}(q_o)\in\tau$ and $\td\g_f:=\td{f}(\td\g_{\,id})$ is the leaf of $\td\G_k$ through $\td{f}(q_o)$.
Since $q_f\in E_k$, but $\td\g_f$ is not contained in $E_k$, $\forall f\in Diff_\G(S,p)$, the immage $\td\Pi_k(\td\g_f):=\g_f$ is an irreducible separatrix of $\G$ through $p$.
Therefore, if we define $\g$ as in (\ref{eq:4}) then $\g$ is a distinguished separatrix of $\G$.
 
We will divide case (2.2) into two subcases:
\begin{itemize}
\item[(2.2.1).] $E_k\setminus\bigcup_j\td{E}_j\ne\emp$.
\item[(2.2.2).] $E_k=\bigcup_j\td{E}_j$.
\end{itemize}

In case (2.2.1), let $\wh{E}$ be a connected component of $E_k\setminus\bigcup_j\td{E}_j$. Let $\bigcup_{i=1}^rD_i$ be decomposition of $\wh{E}$ into irreducible components, $D_i\simeq\p^1$. Note that:
\begin{itemize}
\item[(i).] The graph formed by the divisors $D_i$ is a tree.
\item[(ii).] The intersection matrix $(D_i\,.\,D_j)_{1\le i,j\le r}$ is negative.
\item[(iii).] If $D$ is an irreducible divisor of $E_k$ such that $D\not\sub\wh{E}$ but $D\cap\wh{E}\ne\emp$ then $D=\td{E}_j$ for some $j$. In particular, $D$ is dicritical.
\end{itemize}
In this case, $Sing(\td\G_k)\cap\wh{E}\ne\emp$ and contains a singularity $q$ with a separatrix $\td\g$ not contained in $\wh{E}$.
This is a consequence of Sebastiani's version of Camacho-Sad theorem (see \cite{sb}). In fact, $\td\g$ is not contained in $E_k$, for otherwise it would be contained in some irreducible divisor $D$ of $E_k$ not contained in $\wh{E}$, and $D$ is non dicritical, which contradicts (iii). Therefore, we reduce the problem to case (1).

In case (2.2.2) all irreducible divisors $\td{E}_j$ of $E_k$ are dicritical.
We can assume that $Sing(\td\G_k)=\emp$. In fact, if $q_o\in Sing(\td\G_k)$ then $q_o$ is simple and any of their separatrices cannot be contained in $E_k$, for otherwise some of the divisors $\td{E}_j$ would be non-dicritical. Therefore, we are again in case (1). In particular, we can assume that all divisors $\td{E}_j$ are radial, in the sense that for any $q\in \td{E}_j$ the leaf of $\td\G_k$ through $q$ is transverse to $\td{E}_j$.
Moreover, $m\ge2$ because otherwise $p$ would be a radial singularity of $\G$. 
In particular, we can assume that $\td{E}_1\cap\td{E}_2=\{q_o\}\ne\emp$. Let $\td\g_{\,id}$ be the leaf of $\td\G_k$ through $q_o$.
Note that, for any $f\in Diff_\G(S,p)$ then
\[
q_f:=\td{f}(q_o)=\td{f}(\td{E}_1)\cap\td{f}(\td{E}_2)\in E_k\,\,,
\]
so that $\td\g_f:=\td{f}(\td\g_{\,id})$ is the leaf of $\td\G_k$ through $q_f$.
For each $f\in Diff_\G(S,p)$ the projection $\g_f:=\td\Pi_k(\td\g_f)$ is an irreducible separatrix of $\G$.
Since
\[
A:=\{q_f\,|\,f\in Diff_\G(S,p)\}\sub\bigcup_{i\ne j}\,\td{E}_i\cap\td{E}_j
\]
then $A$ is finite. Therefore, we can construct a distinguished separatrix $\g$ of $\G$ as in (\ref{eq:4}).
\qed

\vskip.1in

Finally, note that lemma \ref{l:21} is a consequence of lemmae \ref{l:24} and \ref{l:25}.

\vskip.1in

\subsection{Proof of lemma \ref{l:22}.}
Since $U$ is a tubular neighborhood of $\Ga$ the map $\Te\in H^2_{Dr}(U)\mapsto\Te|_\Ga\in H^2_{Dr}(\Ga)$ is an isomorphism.
Therefore, it is sufficient to prove that $C_1(N_\Si)|_\Ga=0$.

Recall that the germ of $\fa$ at any $q\in\Ga$ is equivalent to a product of a singular foliation by curves on $(\C^2,0)$ by a regular foliation of dimension $n-2$.
This implies that there exist a local coordinate system around $q$, $z=(x,y)\colon U\to\C^2\times\C^{n-2}$, $x=(x_1,x_2)$, $y=(y_1,...,y_{n-2})$, and a holomorphic vector field $X=P(x)\frac{\pa}{\pa x_1}+Q(x)\frac{\pa}{\pa x_2}$, with an isolated singularity at $0\in\C^2$, such that
\begin{itemize}
\item{} $\fa|_U$ is generated by the $n-1$ commuting vector fields $X$, $Y_1:=\frac{\pa\,}{\pa y_1}$,..., $Y_{n-2}:=\frac{\pa\,\,\,}{\pa y_{n-2}}$.
\end{itemize}
Moreover, the separatrix $\Si$ of $\fa$ along $\Ga$ is induced by a separatrix $\g=(f(x_1,x_2)=0)$ of $X$, such that $X(f)=h.\,f$, where we have assumed $h(0)=0$.
\vskip.1in
It follows that we can find a Leray covering $\U=(U_\a)_{\a\in A}$ of $\Ga$ by open sets with the following properties:
\begin{itemize}
\item[(a).] For each $\a\in A$, there exists a coordinate system $z_\a=(x_\a,y_\a)\colon U_\a\to\C^2\times\C^{n-2}$, where $\Si\cap U_\a=(x_\a=0)$, $x_\a=(x_{\a1},x_{\a2})$ and $y_\a=(y_{\a1},...,y_{\a\,n-2})$.
\item[(b).] For each $\a\in A$, $\fa|_{U_\a}$ is generated by the $n-1$ holomorphic vector fields
\[
X_\a=P(x_\a)\frac{\pa}{\pa x_{\a1}}+Q(x_\a)\frac{\pa}{\pa x_{\a2}}\,,\,Y_{\a j}=\frac{\pa}{\pa y_{\a j}},\,\,j=1,...,n-2\,\,.
\]
\item[(c).] $\Si\cap U_\a$ has the reduced equation $f_\a=0$, where $f_\a=f(x_\a)$. In particular, if we set $h_\a=h(x_\a)$ then
\[
X_\a(f_\a)=h_\a.\,f_a\,,\,Y_{\a j}(f_\a)=0,\,\forall\,1\le j\le n-2\,.
\]
\end{itemize}
Consider the multiplicative cocycle $g=(g_{\a\be})_{U_\a\cap U_\be\ne\emp}$ such that $f_\a=g_{\a\be}.\,f_\be$ on $U_\a\cap U_\be\ne\emp$.
\vskip.1in
\begin{claim}\label{cl:22}
If $U_\a\cap U_\be\ne\emp$ then $g_{\a\be}$ is locally constant on $U_\a\cap U_\be\cap\Ga${\rm :} $dg_{\a\be}|_{U_\a\cap U_\be\cap\Si}\equiv0$. In particular $C_1(N_\Si)=0$.
\end{claim}
{\it Proof.}
Let $U_\a\cap U_\be\cap \Si\ne\emp$. We assert that there exists a $(n-1)\times(n-1)$ matrix $A_{\a\be}$, with entries in $\O(U_\a\cap U_\be)$,
$A_{\a\be}=\left(a_{\a\be}^{ij}\right)_{0\le i,j\le n-2}$, such that
\begin{equation}
\left\{
\begin{matrix}\label{eq:5}
X_\a=a_{\a\be}^{00}.\,X_\be+\sum_{j=1}^{n-2}a_{\a\be}^{0j}.\,Y_{\be j}&\,\\
Y_{\a i}=a_{\a\be}^{i0}.\,X_\be+\sum_{j=1}^{n-2}a_{\a\be}^{ij}.\,Y_{\be j}&,\,1\le i\le n-2\\
\end{matrix}
\right.
\end{equation}

In fact, since $\fa|_{U_\a\cap U_\be}$ is generated by both systems $\left<X_\a,Y_{\a i}\,|\,1\le i\le n-2\right>$ and $\left<X_\be,Y_{\be i}\,|\,1\le i\le n-2\right>$, we can find a matrix $A_{\a\be}$ with entries in $\O(U_\a\cap U_\be\setminus\Si)$ as in (\ref{eq:5}). But since $cod(\Si)=2$ the entries of $A_{\a\be}$ can be extended to $U_\a\cap U_\be$ by Hartog's theorem.

Now, from (c) we get
\[
0=Y_{\a i}(f_a)=Y_{\a i}(g_{\a\be}.\,f_\be)=Y_{\a i}(g_{\a\be}).\,f_\be+g_{\a\be}.\,Y_{\a i}(f_\be)
\]
and from (c) and (\ref{eq:5})
\[
Y_{\a i}(f_\be)=a_{\a\be}^{i0}.\,X_\be(f_\be)+\sum_{j=1}^{n-2}a_{\a\be}^{ij}.\,Y_{\be j}(f_\be)=a_{\a\be}^{i0}.\,h_\be.\,f_\be\,\implies
\]
\[
\left(Y_{\a i}(g_{\a\be})+a_{\a\be}^{i0}.\,h_\be\right)\,f_\be=0\,\implies\,Y_{\a i}(g_{\a\be})=-a_{\a\be}^{i0}.\,h_\be\,.
\]
Now, $h_\be|_{U_\a\cap U_\be\cap\Si}=h(0)=0$ and so
\[
Y_{\a i}(g_{\a\be})|_{U_\a\cap U_\be\cap\Si}=\frac{\pa g_{\a\be}}{\pa y_{\a i}}(0,y_\a)=0\,,\,1\le i\le n-2\,\,\implies\,\,dg_{\a\be}|_{U_\a\cap U_\be\cap\Si}=0\,.
\]
This finishes the proof of lemma \ref{l:22}.
\qed

\vskip.1in

\subsection{Proof of lemma \ref{l:23}.}\label{ss:23}
The case in which $\Si$ is smooth was proved in \cite{lc2}. 
Here we give a more general proof (suggested by J. V. Pereira).
Let us consider first the case $n=3$: $M=\p^3$.
In this case, $\Ga$ is a compact algebraic curve so that $H^2_{Dr}(\Ga)\simeq\R$ and the map
\[
\Te\in H^2_{Dr}(\Ga)\mapsto\int_\Ga\Te\in\R
\]
is an isomorphism.
In fact, we will prove that
\[
\int_\Ga C_1(N_\Si)\in\N\,\,\implies\,\,C_1(N_\Si)\ne0\,\,.
\]
We will see that $\int_\Ga C_1(N_\Si)$ represents the intersection number of a small deformation $\Ga_t$ of $\Ga$ with $\Si$.

Let $\X(\p^3)$ be the vector space of holomorphic vector fields on $\p^3$: $dim(\X(\p^3))=15$.
Given $Z\in\X(\p^3)$ we will denote by $(t,q)\in\C\times\p^3\mapsto Z_t(q)\in\p^3$ its flow and $\Ga_t:=Z_t(\Ga)$.
Let $U$ be a tubular neighborhood $U$ of $\Ga$ with $U\sub\bigcup_\a U_\a$.
\begin{rem}\label{r:26}
{\rm There exist $Z\in\X(\p^3)$ and $\ep>0$ with the following properties:
\begin{itemize}
\item[(a).] If $t\in D_\ep\sub\C$ then $\Ga_t\sub U$, where $D_\ep=\{t\,|\,|t|<\ep\}$.
\item[(b).] If $t\in D_\ep^*:=D_\ep\setminus\{0\}$ then $\Ga_t\cap\Ga=\emp$.
\item[(c).] The set $B:=\{t\in D_\ep\,|\,\Ga_t\,$ is not transverse to $\Si\}$ is discrete in $D^*$.
\item[(d).] There exists $t_o\in D^*\setminus B$ such that $\Ga_{t_o}\cap\Si\ne\emp$.
\end{itemize}}
\end{rem}
We leave the proof of remark \ref{r:26} for the reader.
Let us finish the proof of lemma \ref{l:23} in the case of $\p^3$.

The idea is to prove that, if $t\in D_\ep^*\setminus B$ then $\int_\Ga C_1(N_\Si)=\#(\Ga_t\cap\Si)$, the intersection number of $\Ga_t$ with $\Si$.
By (d) of remark \ref{r:26} $\#(\Ga_t\cap\Si)>0$ and so $C_1(N_\Si)\ne0$.

\vskip.1in

First of all, note that
\[
\Ga_t\cap\Si=Z_t(\Ga\cap Z_{-t}(\Si))\,\implies\,\#\,[\Ga_t\cap\Si]=\#\,[\Ga\cap Z_{-t}(\Si)]\,.
\]
On the other hand, $Z_{-t}(\Si)$ can be defined in the covering $\U_t:=(Z_{-t}(U_\a))_{\a\in A}$ by the divisor $(f_\a\circ Z_t)_{\a\in A}$, with associated cocycle
$g_t:=(g_{\a\be}\circ Z_t)_{U_\a\cap U_\be\ne\emp}$.
Since $t\in D^*_\ep\setminus B$, $\Ga$ is transverse to $Z_{-t}(\Si)$ and so $\Ga\cap Z_{-t}(\Si)$ is finite and is defined by the divisor
$\left(f_\a\circ Z_t|_{\Si\cap Z_{-t}(U_\a)}\right)_{\a\in A}$ with associated cocycle
$g_t|_\Ga=\left(g_{\a\be}\circ Z_t|_{\Si\cap Z_{-t}(U_\a\cap U_\be)}\right)_{U_\a\cap U_\be\ne\emp}$. This divisor can be interpreted as a holomorphic section of the line bundle induced by $g_t|_\Ga$ on $Picc(\Ga)$.
In particular, if $C_1(g_t|_\Ga)$ is its first Chern class then its degree is done by
\[
\int_\Ga C_1(g_t|_\Si)=\#\,[\Ga\cap Z_{-t}(\Si)]=\#\,[\Ga_t\cap \Si]\,\,.
\]
Since the map
\[
t\in D_\ep\mapsto \int_\Ga C_1(g_t|_\Ga)
\]
is continuous and constant in $D_\ep^*\setminus B$, we get $\int_\Ga C_1(g_0|_\Si)>0$ $\implies$ $C_1(g_0|_\Si)=C_1(N_\Si)\ne0$.
This finishes the proof of lemma \ref{l:23} in the case of $\p^3$.

The case of $\p^n$, $n\ge4$, can be reduced to the previous by taking sections by generic 3-planes linearly embedded in $\p^n$.
We leave the details to the reader.
\qed

\subsection{Proof of corollary \ref{c:21}.}
Recall that $X(f)=h.\,f$, where $X$ represents the normal type $\G$ of $\fa$ along $\Ga$ and $f\in\O_2$ is reduced.
By lemma \ref{l:23} we have $h(0)\ne0$.
Let $f_\mu$ and $X_\nu$ be the first non-zero jets of $f$ and $X$ at $0\in\C^2$, respectively.
Then
\[
X(f)=h.\,f\,\,\implies\,\,X_\nu(f_\mu)=h(0).\,f_\mu\,\,\implies\,\,\nu=1
\]
and $X_\nu=X_1$ is not nilpotent; has at least one non-zero eigenvalue.
On the other hand, we have seen that $\Ga$ is a Kupka component of $\fa$ if, and only if, $tr(X_1)\ne0$. If $tr(X_1)=0$ and $X_1$ has a non-zero eigenvalue, then we can assume that
$X_1=\la\,\left(x_1\frac{\pa}{\pa x_1}-x_2\frac{\pa}{\pa x_2}\right)$, $\la\ne0$. In this case, $X$ has exactly two separatrices through $0\in\C^2$ which are smooth and tangent to $x_1=0$ and $x_2=0$. We can assume that these separatrices have equations $f_1(x_1,x_2)=x_1+h.o.t$ and $f_2(x_1,x_2)=x_2+h.o.t$.
Consider the separatrix $\g=(f_1.\,f_2=0)$ of $X$.
Note that $f(\g)=\g$, $\forall f\in Diff_\G(\C^2,0)$. 
By lemma \ref{l:24} $\g$ generates a separatrix $\Si$ of $\fa$ along $\Ga$.
However $X(f_1.\,f_2)=h.\,f_1.\,f_2$ where $h(0)=0$, because $X_1(x_1.\,x_2)=0$.
Therefore, we must have $tr(X_1)\ne0$ and $\Ga$ is a Kupka component of $\fa$.
\qed 

\section{Corollaries \ref{c:2} and \ref{c:3}}\label{ss:3}

\subsection{Proof of corollary \ref{c:2}.}
A codimension one foliation $\G$ on $\p^n$ of degree zero has a rational first integral of degree one. It is defined in some coordinate system $(x_1,...,x_{n+1})\in\C^{n+1}$ by a the form
$\om=x_1\,dx_2-x_2\,dx_1$. In particular, $\Pi^{-1}(Sing(\G))=(x_1=x_2=0)$, which is a l.t.p component.

\vskip.1in

Conversely, let $\fa$ be a codimension one foliation on $\p^n$, $n\ge3$. It is known that $\fa$ has at least one irreducible component of codimension two \cite{ln}. Assume that all components of $Sing(\fa)$ are l.t.p.
Let $\Om$ be a 1-form on $\C^{n+1}$ that represents $\fa$ in homogeneous coordinates: $\fa_\Om=\Pi^*(\fa)$. Then
\begin{itemize}
\item[(a).] $i_R\Om=0$, where $R$ is the radial vector field on $\C^{n+1}$.
\item[(b).] The coefficients of $\Om$ are homogeneous of degree $d+1$, where $d=deg(\fa)$.
\item[(c).] $i_Rd\Om=(d+2)\,\Om$ (see \cite{lc}). In particular, $Sing(d\Om)\sub Sing(\Om)$.
\end{itemize}

\begin{claim}\label{cl:31}
Let $q\in Sing(\fa)$ and $p\in\Pi^{-1}(q)\sub\C^{n+1}\setminus\{0\}$. Then $d\Om_p\ne0$.
In particular, $Sing(d\Om)=\emp$ and $deg(\fa)=0$.
\end{claim}

{\it Proof.}
Let $\om$ be a holomorphic 1-form that represents $\fa$ in a neighborhood of $q$.
The hypothesis and theorem \ref{t:1} imply that $d\om(p)\ne0$.

On the other hand, $\Pi^*(\om)$ represents $\fa_\Om$ in a neighborhood, say $U$, of $p$.
It follows that $\Pi^*(\om)=\var.\,\Om$ on $U$, where $\var\in\O^*(U)$.
Therefore,
\[
\Pi^*(d\om)=d\,\Pi^*(\om)=d\var\wedge\Om+\var.\,d\Om\,\,\implies
\]
\[
\forall\,u,v\in T_p\C^{n+1}\,\,\text{we get}\,\,\var(p).\,d\Om_p(u,v)=\Pi^*(d\om)_p\,(u,v)=d\om_q\,(d\Pi(p).u,d\Pi(p).v)\,.
\]
Since $\Pi$ is a submersion, it follows that $d\Om_p\ne0$.
Therefore, the coefficients of $\Om$ must be of degree one and $\fa$ has degre zero, as asserted in corollary \ref{c:2}.
\qed

\vskip.1in

\subsection{Proof of corollary \ref{c:3}.}
The idea is to use corollary \ref{c:2}.
Assume that there exists an integrable 2-form $\eta$ on $\C^n$, $n\ge4$, with homogeneous coefficients of degree $d\ge1$ and such that $Sing(\eta)=\{0\}$.
Denote by $\fa_\eta$ the holomorphic codimension two foliation of $\C^n$ generated by $\eta$.
By assumption $Sing(\fa_\eta)=\{0\}$.
Note also that the codimension two distribution of $\C^n\setminus\{0\}$ tangent to $\fa_\eta$ is given by
\[
ker(\eta)(p)=\{v\in T_p\C^n\,|\,i_v\,\eta(p)=0\}\,,\,\forall p\ne0\,,
\]
where $i_v$ denotes the interior product.
The fact that $ker(\eta)$ has codimension two is equivalent to
\begin{equation}\label{eq:e2}
\eta\wedge\eta=0\,\,.
\end{equation}

Let $\om=i_R\eta$, where $R=\sum_{j=1}^nz_j\frac{\pa}{\pa z_j}$ is the radial vector field on $\C^n$.
We have two possibilities: either $\om\equiv0$, or $\om\not\equiv0$.

In the first case, $\eta$ generates a codimension two foliation on $\p^{n-1}$: there exists a codimension two foliation $\fa$ on $\p^{n-1}$ such that $\Pi^*(\fa)=\fa_\eta$, where $\Pi\colon\C^n\setminus\{0\}\to\p^{n-1}$ denotes the canonical projection. However, any codimension two foliation on $\p^{n-1}$, $n\ge4$, has at least one singularity: if
$q\in Sing(\fa)$ then the line $\ov{\Pi^{-1}(q)}\sub\C^n$ is contained in the singular set of $\eta$.

In the second case $\om$ is a 1-form on $\C^n$ with homogeneous coefficients of degree $d+1$.

\begin{lemma}\label{int}
The form $\om$ is integrable: $\om\wedge d\om=0$.
\end{lemma}

{\it Proof.}
The following is equivalent to the integrability of the distribution $ker(\eta)$:
\begin{itemize}
\item[(I).] for any $p\in \C^n\setminus\{0\}$ there exists a germ coordinate system $(x,y)\colon(\C^n,p)\to(\C^2,0)\times(\C^{n-2},0)$, with $x=(x_1,x_2)$, such that
$\eta_p=\var(x,y)\,dx_1\wedge dx_2$, where $\eta_p$ is the germ of $\eta$ at $p$ and $\var\in\O_p^*$.
\end{itemize}
Since the coefficients of $\eta$ are homogeneous of degree $d$ we have $L_R\eta=(d+2)\,\eta$, where $L_R$ denotes the Lie derivative in the direction of $R$.
From this we get
\[
(d+2)\,\eta=L_R\eta=i_Rd\eta+d\,i_R\eta=i_Rd\eta+d\om\,\,\implies
\]
\[
\om\wedge d\om=i_R\eta\wedge d\om=(d+2)\,i_R\eta\wedge \eta-i_R\eta\wedge i_Rd\eta\,\,.
\]
Now, from (\ref{eq:e2}) we get
\[
0=i_R(\eta\wedge\eta)=2\,i_R\eta\wedge\eta=2\,\om\wedge\eta\,\implies\,\om\wedge d\om=-i_R\eta\wedge i_Rd\eta\,.
\]
If we consider a coordinate system as in (I) we have $\eta_p=\var\,dx_1\wedge dx_2$ and $d\eta_p=d\var\wedge dx_1\wedge dx_2$ and this 
implies that $i_R\eta_p\wedge i_Rd\eta_p=0$, as the reader can check.
\qed
\vskip.1in
Write $\om=\phi.\,\om_1$, where $\phi$ is homogeneous and $cod(Sing(\om_1))\ge2$.
\begin{rem}\label{r:31.}
{\rm 
Note that:
\begin{itemize}
\item[(a).] $\om_1\wedge\eta=0$. This is a consequence of $\om\wedge\eta=0$.
\item[(b).] $\om_1\wedge d\om_1=0$ and $i_R\om_1=0$. This is a consequence of $\om\wedge d\om=0$ and $i_R\om=0$.
\end{itemize}}
\end{rem}
Denote by $\fa_{\om_1}$ the foliation generated by $\om_1$.
It follows from (b) of remark \ref{r:31.} that there exists a codimension one foliation $\fa$ on $\p^{n-1}$ such that $\Pi^*(\fa)=\fa_{\om_1}$.

\begin{lemma}\label{l:32}
All irreducible components of $Sing(\fa)$ are l.t.p.
\end{lemma}

{\it Proof.}
Fix $q\in Sing(\fa)$ and $p\in\C^n\setminus\{0\}$ with $\Pi(p)=q$. Note that $p\in Sing(\fa_{\om_1})$, the foliation generated by $\om_1$.
Let $(x,y)\colon(\C^n,p)\to(\C^2,0)\times(\C^{n-2},0)$ be as in (I), so that $\eta=\var.\,dx_1\wedge dx_2$, $\var\in\O_p^*$.
It follows from $\om_1\wedge\eta=0$ that in these coodinates we have $\om_1=A(x,y)\,dx_1+B(x,y)\,dx_2$ and from $\om_1\wedge d\om_1=0$ that
\[
(A\,dB-B\,dA)\wedge dx_1\wedge dx_2=0\,\,\implies\,\,\om_1=h(x,y).\,\left(C(x_1,x_2)\,dx_1+D(x_1,x_2)\,dx_2\right)\,.
\]
Since $cod(Sing(\om_1))\ge2$ we get $h\in\O_p^*$ and the germ of $Sing(\om_1)$ at $p$ is defined by $(x_1=x_2=0)$.
Moreover, the germ of $\fa_{\om_1}$ at $p$ is defined by the form $C(x_1,x_2)\,dx_1+D(x_1,x_2)\,dx_2$ and so $\fa_{\om_1}$ is a transversely product at $p$.
Since $p\in\Pi^{-1}(q)$ and $\Pi$ is a submersion at $p$, $\fa$ is a transversely product at $q$.
\qed

\vskip.1in

Corollary \ref{c:2} implies that $\om_1$ has a linear rational first integral that we can assume to be $x_2/x_1$, so that $\om_1=x_1\,dx_2-x_2\,dx_1=i_R(dx_1\wedge dx_2)$.
Let $\eta=\sum_{i<j}\eta_{ij}\,dx_i\wedge dx_j$, where $\eta_{ij}$ is homogeneous of degree $d$, $\forall\,i<j$.
From $\om_1\wedge\eta=0$ we get $\eta_{ij}=0$, $\forall\,j>i\ge3$. Therefore, we can write
$\eta=dx_1\wedge\a+dx_2\wedge\be+\g\,dx_1\wedge dx_2$, where $\a=\sum_{j\ge3}\eta_{1j}\,dx_j$, $\be=\sum_{j\ge3}\eta_{2j}\,dx_j$ and $\g=\eta_{12}$.
Hence,
\[
0=\om_1\wedge\eta=(x_1\,dx_2-x_2\,dx_1)\wedge(\a\wedge dx_1+\be\wedge dx_2+\g\,dx_1\wedge dx_2)\,\implies
\]
\[
(x_1\,\a+x_2\,\be)\wedge dx_1\wedge dx_2=0\,\,\implies\,\,x_1\,\a=-x_2\,\be\,\,\implies
\]
there exists 1-form $\mu$ with homogeneous coefficients of degree $d-1$ such that $\a=-x_2\,\mu$ and $\be=x_1\,\mu$.
In particular, we get
\[
\eta=\om_1\wedge\mu+\g\,dx_1\wedge dx_2=(x_1\,dx_2-x_2\,dx_1)\wedge\mu+\g\,dx_1\wedge dx_2\,\implies\,
\]
\[
Sing(\eta)\sup(x_1=x_2=\g=0)\,\implies
\]
$d=0$ and $\g$ is a constant, for otherwise $cod(Sing(\eta))\le3$ and $Sing(\eta)\supsetneq\{0\}$.
This finishes the proof of corollary \ref{c:3}.
\qed

\vskip.2in

\bibliographystyle{amsalpha}

\vskip.3in

{\sc A. Lins Neto}

{\sl Instituto de Matem\'atica Pura e Aplicada}

{\sl Estrada Dona Castorina, 110}

{\sl Horto, Rio de Janeiro, Brasil}

{\tt E-Mail: alcides@impa.br}


\end{document}